\documentclass{amsart}

\newtheorem{theorem}{Theorem}
\newcommand{\bt}{\begin{theorem}}
\newcommand{\et}{\end{theorem}}
\newtheorem{lemma}{Lemma}
\newcommand{\bl}{\begin{lemma}}
\newcommand{\el}{\end{lemma}}
\newtheorem{corollary}{Corollary}
\newcommand{\bc}{\begin{corollary}}
\newcommand{\ec}{\end{corollary}}
\newcommand{\beq}{\begin{equation}}
\newcommand{\eeq}{\end{equation}}
\newcommand{\benum}{\begin{enumerate}}
\newcommand{\eenum}{\end{enumerate}}
\newcommand{\F}{\ensuremath{\mathbf F }}
\newcommand{\N}{\ensuremath{ \mathbf N }}
\newcommand{\Z}{\ensuremath{\mathbf Z}}

\title[Growth polynomials for additive $(h,k)$-tuples]
{Growth polynomials for additive quadruples and  $(h,k)$-tuples} 

\author{Melvyn B. Nathanson}
\address{Lehman College (CUNY),Bronx, New York 10468}
\email{melvyn.nathanson@lehman.cuny.edu}

\subjclass[2010]{ .}

\thanks{Supported in part by a grant from the PSC-CUNY Research Award Program.}

\begin{document}

\maketitle

\begin{abstract}
Consider the interval of integers $I_{m,n} = \{ m, m+1, m+2,\ldots, m+n-1 \}$.
For fixed integers $h,k,m$, and $c$, let $\Phi_{h,k,m}^{(c)}(n)$ denote the number of solutions
of the equation  $(a_1+\cdots + a_h )- (a_{h+1} + \cdots + a_{h+k})=c$
with $a_i \in I_{m,n}$ for all $i=1,\ldots, h+k$.   
This is a polynomial in $n$ for all sufficiently large $n$, 
and the growth polynomial is constructed explicitly.  
\end{abstract}

\section{Linear relations in the integers}

Let $G$ be an additive abelian group, and let $A$ be a nonempty finite subset of $G$.
Let $A^h$ denote the set of $h$-tuples of elements of $A$.
An \emph{additive quadruple} in $A$ is a 4-tuple $(a_1,a_2,a_3,a_4) \in A^4$ 
such that $a_1+a_2=a_3+a_4$.   
These quadruples are important in additive number theory 
because a set $A$ with a small sumset $2A = \{a_1+a_2:a_1,a_2 \in A\}$ 
must have many additive 
quadruples, and a set with few additive quadruples must have 
a large sumset  (cf. Nathanson~\cite{nath96bb} and Tao-Vu~\cite{tao-vu06}).  
For any positive integer $h$,  define $\Psi_h(A)$ as the number of $2h$-tuples 
$(a_1,\ldots, a_{2h}) \in A^{2h}$ such that 
\beq    \label{quad:DE=}
a_1 + a_2 + \cdots +  a_{h} = a_{h+1} + a_{h+2} + \cdots + a_{2h}. 
\eeq
The following problem comes from counting the number of additive quadruples 
and the number of solutions of  more general linear  relations 
in the additive group of integers.

Let  \N\ and \Z\ denote the sets of positive integers and integers, respectively.
For $n \in \N$ and $m\in \Z$, we define the intervals of  integers 
\[
I_n = \{0,1,2,\ldots, n-1\}
\]
and
\[
I_{m,n} = \{ m, m+1, m+2,\ldots, m+n-1 \}.
\]
For $h,k \in \N$ and $c \in \Z$, let $\Phi_{h,k}^{(c)}(n)$ denote the number 
of $(h+k)$-tuples 
\[
(a_1,a_2,\ldots, a_h, a_{h+1}, a_{h+2},\ldots, a_{h + k} ) \in I_{n}^{h +k}
\]
such that 
\beq    \label{quad:DEc}
a_1 + a_2 + \cdots +  a_{h} - a_{h+1} - a_{h+2} - \cdots - a_{h + k} = c. 
\eeq
For $h,k \in \N$ and $m,c \in \Z$, let $\Phi_{h,k,m}^{(c)}(n)$ denote the number 
of $(h+k)$-tuples in $I_{m,n}^{h +k}$ that satisfy equation~\eqref{quad:DEc}.  
Thus, $\Phi_{h,k}^{(c)}(n) = \Phi_{h,k,0}^{(c)}(n)$.

The arithmetic function 
\[
\Psi_h(n) = \Psi_h(I_n) = \Phi_{h,h}^{(0)}(n) 
\]
counts the number of solutions of  equation~\eqref{quad:DE=} 
with $a_i \in I_n$ for $i=1,\ldots, 2h$.

Note that $\Phi_{h,k,m}^{(c)}(n) > 0$ if and only if 
\[
-kn + k + (h-k)m \leq c \leq hn - h + (h-k)m.
\]
Thus, for fixed $h,k,m$, and $c$ we have $\Phi_{h,k,m}^{(c)}(n) > 0$ 
for all sufficiently large $n$.  

We shall prove that, for all $h, k \in \N$ and $m,c\in \Z$, the arithmetic function 
$\Phi^{(c)}_{h,k,m}(n)$ is a polynomial in $n$ of degree 
$h+k - 1$ for all sufficiently large $n$.
In particular, 
\[
\Psi_h(n) = \Phi^{(0)}_{h,h}(n)
\]
is a polynomial in $n$ of degree 
$2h - 1$ for all sufficiently large $n$.

If $h =k$ and if $(a_1,\ldots, a_{2h}) \in \Z^{2h}$ is a solution
of equation~\eqref{quad:DEc}, then the translate 
$(a_1+t,\ldots, a_{2h}+t) \in \Z^{2h}$ 
is also a solution of equation~\eqref{quad:DEc} for all $t \in \Z$.  
This implies that $\Phi^{(c)}_{h,h, m}(n)= \Phi^{(c)}_{h,h}(n)$ for all $m,c \in \Z$.

\section{Polynomial growth of additive quadruples}
We begin by explicitly computing the  polynomials  $\Phi_{h,1}^{(0)}(n)$ 
and $\Phi_{2,2}^{(0)}(n)$.  
Recall two standard facts about binomial coefficients.  

If $r_h(x)$ denotes the number of $h$-tuples $(a_1,\ldots, a_h)$ 
of  nonnegative 
integers such that $a_1 + \cdots + a_h = x$,  then 
\[
r_h(x) = \binom{x + h-1}{h -1}
\]
and
\[
\sum_{x=0}^{n-1} \binom{x + h-1}{h -1} = \binom{n + h-1}{h }.
\]
We call $r_h(x)$ the $h$-fold \emph{representation function}.  

Let $d$ be a nonnegative integer. 
If $a,b,c \in \Z$ and $(a,b)\neq (0,0)$, then the binomial coefficient 
\begin{align*}
\binom{ax+by+c}{d} 
& = \frac{1}{d!}(ax+by+c)(ax+by+c-1)\cdots (ax+by+c-d+1) \\
& = \frac{(ax+by)^d}{d!} + \cdots + \binom{c}{d}.
\end{align*}
is a polynomial in $x$ and $y$ of total degree $d$.  
If $b=0$, then this is a polynomial in $x$ of degree $d$.

\bt
Let $h \in \N$.  For all $n \in \N$,
\[
\Phi^{(0)}_{h,1}(n)  = \binom{n+h-1}{h}  = \frac{1}{h!}n^h + \cdots.
\]
is a polynomial in $n$ of degree $h$ with leading coefficient $1/h!$.
\et

\begin{proof}
Counting the number of $(h+1)$-tuples in 
$(a_1,a_2,\ldots, a_h, a_{h+1}) \in I_n^{h+1}$ such that 
\[
a_1 + a_2 + \cdots +  a_{h} = a_{h+1}
\]
we obtain 
\[
\Phi^{(0)}_{h,1}(n)  = \sum_{a_{h+1}=0}^{n-1} r_h(a_{h+1}) 
= \sum_{a_{h+1}=0}^{n-1}  \binom{a_{h+1}+h-1}{h-1} = \binom{n+h-1}{h} 
\]
and so $ \Phi_{h,1}^{(0)}(n)$ is a polynomial in $n$ of degree $h$ 
with leading coefficient $1/h!$.
This completes the proof. 
\end{proof}

\bt[Additive quadruples]   \label{quad:theorem:(2,2)poly}
For all $n \in \N$,
\[
\Psi_2(n)  = \frac{2}{3}n^3 + \frac{1}{3}n
\]
is a polynomial in $n$ of degree 3.
\et

The integer $(2n^3+n)/3$ is the $n$th \emph{octahedral number}.  

\begin{proof}
Here are two proofs.  

The function $\Psi_2(n)$ counts the number of quadruples 
$(a_1, a_2, a_3, a_4)$ such that 
\[
a_1 + a_2 = a_3 + a_4
\]
and $a_i \in I_n = \{0,1,2,\ldots, n-1\}$ for $i=1,2,3,4$.  
Because $a_1+a_2=x$ if and only if $(n-1-a_1) + (n-1-a_2) = 2n-2-x$, 
and because $r_2(x) = x+1$ for $x = 0,1,\ldots, n-1$, it follows that 
\begin{align*}
\Psi_2(n)  
& = r_2(n-1)^2  + 2 \sum_{x=0}^{n-2} r_2(x)^2 \\
& =  n^2 + 2 \sum_{x=0}^{n-2} (x+1)^2\\
& =  n^2  + \frac{1}{3} (n-1)n(2n-1)  \\
& = \frac{2}{3}n^3 + \frac{1}{3}n.
\end{align*}

The second proof is by induction on $n$.  
We have $\Psi_2(1)  = 1$ because 0 + 0 = 0 + 0.
Suppose that the Theorem holds for some positive integer $n$.
For $j  \in \{ 0,1,2,3,4 \}$, let $N(j)$ denote the number of quadruples 
$(a_1,a_2,a_3,a_4)\in I_{n+1}^4$ such that 
$a_1+a_2=a_3+a_4$ and $a_i = n$ for exactly $j$ integers $i \in \{1,2,3,4\}$.
Then $\Psi_2(n+1) = \sum_{j=0}^4 N(j)$.  
We have $N(0) = \Psi_2(n)$, $N(3)=0$, and $N(4)=1$.

The only quadruples with $j=2$ are of the form 
$(a,n,a,n)$, $(a,n,n,a)$, $(n,a,a,n)$, and $(n,a,n,a)$ 
with $a \in I_n$, and so $N(2)=4n$.

Quadruples with $j=1$ must be of the form
$ (a_1,a_2,a_3,n), (a_1,a_2,n,a_4), (a_1,n,a_3,a_4)$, or $(n,a_2,a_3,a_4)$, 
with $a_i \in I_n$, and there is the same number of quadruples of each form.    
Consider quadruples of the form $ (a_1,a_2,a_3,n)$.
Then $a_1+a_2=a_3+n$, and $0 \leq a_3 \leq n-2$ because $a_1+a_2 \leq 2n-2$.  
For each $a_3$ we must have $(a_1,a_2) = (x,a_3+n-x)$ with $0 \leq x \leq n-1$ and 
$0 \leq a_3+n-x \leq n-1$.  
Equivalently, $a_3+1 \leq x \leq n-1$, and this interval contains 
exactly $n-1 - a_3$ integers $x$.
It follows that the number of quadruples of the form 
$ (a_1,a_2,a_3,n)$ is 
\[
\sum_{a_3=0}^{n-2} (n-1 - a_3) = \frac{n(n-1)}{2}
\]
and so 
\[
N(1) = 4\frac{n(n-1)}{2} = 2n^2 - 2n.
\]
Thus, 
\begin{align*}
\Phi_{2,2}^{(0)}(n+1) & = \Phi_{2,2}^{(0)}(n) + \sum_{j=1}^4 N(j) \\ 
& = \left( \frac{2}{3}n^3  + \frac{1}{3}n  \right)+ \left( 2n^2 - 2n \right) + 4n + 0 + 1 \\ 
& = \frac{2}{3} (n+1)^3  + \frac{1}{3}(n+1).
\end{align*}
This completes the second proof.  
\end{proof}

\section{Explicit growth polynomials}

In this section we prove that, for fixed $h, k, m$, and $c$,  
the function $\Phi^{(c)}_{h,k,m}(n)$ is, for sufficiently large integers $n$,  
a polynomial in $n$ of degree 
exactly $h+k -1$.  
The first step is to show that  $\Phi^{(c)}_{h,k,m}(n)$ 
has order of magnitude $n^{h+k -1}$.  

We denote by $[x]$ the integer part of the real number $x$.

\bl             \label{quad:lemma:polygrowth}
For $h,k \in \N$ and $m,c \in \Z$, there exists a number $\theta = \theta(h,k,m,c) > 0$  
such that 
\[
\theta n^{h+k - 1} \leq \Phi^{(c)}_{h,k,m}(n) \leq  n^{h+k - 1} 
\]
for all sufficiently large integers $n$.
\el

\begin{proof}
Because the equations
\[
 \sum_{i=1}^h a_i - \sum_{i=h+1}^{h+k} a_i = c
\qquad
\text{and}
\qquad
 \sum_{i=h+1}^{h+k} a_i  - \sum_{i=1}^h a_i = -c
\]
have the same solutions in $I_{m,n}$, it follows that
\[
\Phi^{(c)}_{h,k,m}(n)  = \Phi^{(-c)}_{k,h,m}(n).
\]
Thus, we can assume without loss of generality that $h \geq k$.

For every  $(h+k - 1)$-tuple $(a_1,\ldots, a_{h+k - 1}) \in I_{m,n}^{h+k -1}$,  
there exists at most one integer $a_{h+k} \in I_{m,n}$ that satisfies~\eqref{quad:DEc}
and so 
\[
\Phi^{(c)}_{h,k,m}(n) \leq  \left| I_{m,n}\right|^{h+k - 1} = n^{h+k - 1}.
\]

Let 
\[
0 < \varepsilon < \frac{1}{3}
\]
and choose a positive integer $N$ such that 
\beq        \label{quad:epsilon}
\varepsilon N > \max( |c|, (h-k)m).
\eeq
Let $\alpha, \beta, \gamma$, and $\delta$ be real numbers such that  
\beq    \label{quad:alpha1}
0 \leq \alpha < \beta  < 1\qquad \text{and} \qquad  0 \leq \gamma < \delta < 1.
\eeq 
Suppose that $(a_1,\ldots, a_{h+k - 1}) \in I_{m,n}^{h+k -1}$ and 
\beq        \label{quad:mn1}
m + \alpha n \leq a_i \leq m + \beta n  \qquad \text{for $i=1,\ldots, h$}
\eeq
and
\beq    \label{quad:mn2}
m + \gamma n \leq a_i \leq m + \delta n  \qquad \text{for $i= h+1,\ldots, h+k-1$.}
\eeq
We define
\[
a_{h+k} = \sum_{i=1}^h a_i - \sum_{i=h+1}^{h+k-1} a_i - c.
\]
If $n \geq N$ and if
\beq    \label{quad:alpha2}
h\beta - (k -1) \gamma + 2\varepsilon \leq 1
\eeq
then
\begin{align*}
a_{h+k} & = \sum_{i=1}^h a_i - \sum_{i=h+1}^{h+k-1} a_i - c \\
& \leq h(m + \beta n )  - (k -1)(m + \gamma n ) -c  \\
& = m + (h\beta - (k -1) \gamma  ) n  + (h-k)m -c\\
& < m + (h\beta - (k -1) \gamma ) n + 2\varepsilon N\\
& \leq m + (h\beta - (k -1) \gamma + 2\varepsilon ) n \\
& \leq m+n.
\end{align*}
Similarly, if $n \geq N$ and if
\beq    \label{quad:alpha3}
h\alpha - (k -1) \delta  - \varepsilon \geq 0
\eeq
 then
\begin{align*}
a_{h+k} 
& \geq h(m + \alpha n )  - (k -1)(m + \delta n ) -c  \\
& =  m + (h\alpha - (k -1) \delta  ) n  + (h-k)m -c\\
& \geq m + (h\alpha - (k -1) \delta  ) n -c\\
& >   m + (h\alpha - (k -1) \delta - \varepsilon ) n \\
& \geq m.
\end{align*}
Thus, inequalities~\eqref{quad:epsilon}--\eqref{quad:alpha3} 
imply that $m \leq a_{h+k} < m+n$, that is, $a_{h+k} \in I_{m,n}$.

We shall construct real numbers $\alpha, \beta, \gamma$, and $\delta$ 
that satisfy inequalities~~\eqref{quad:alpha1},~\eqref{quad:alpha2}, and~\eqref{quad:alpha3}.

There are two cases.  The first case is $k = 1$.  
Because $0< \varepsilon  < 1/3$, we have $\varepsilon < 1 - 2\varepsilon$ 
and  there exist numbers $\alpha, \beta \in (0,1)$ such that 
\[
\frac{\varepsilon}{h} \leq \alpha < \beta \leq \frac{1 - 2\varepsilon}{h}.
\]
Let $\gamma$ and $\delta$ be any numbers such that $0 \leq \gamma < \delta < 1$.
Inequalities~\eqref{quad:alpha1},~\eqref{quad:alpha2}, 
and~\eqref{quad:alpha3} are satisfied.  

In the second case, we have  $k \geq 2$.  
Because $0< \varepsilon  < 1/3$, there exist $\gamma, \delta \in (0,1)$ such that 
\[
0 < \delta - \gamma  < \frac{1-3\varepsilon}{k - 1}.
\]
It follows that 
\[
0 < (k - 1)\delta + \varepsilon < (k - 1)\gamma +  1 - 2 \varepsilon < k \leq h
\]
and so there exist numbers $\alpha, \beta \in (0,1)$ such that 
\[
0 < \frac{ (k - 1)\delta  + \varepsilon }{h} < \alpha < \beta 
< \frac{ (k - 1)\gamma + 1 - 2\varepsilon }{h} < 1.
\]
Thus, $\alpha,\beta,\gamma$, and $\delta$ satisfy
conditions~\eqref{quad:alpha1},~\eqref{quad:alpha2}, and~\eqref{quad:alpha3}. 

If the $(h + k - 1)$-tuple $(a_1,\ldots, a_{h+k-1})$ 
satisfies conditions~\eqref{quad:mn1} and~\eqref{quad:mn2}, 
then there exists a unique integer $a_{h+k} \in I_{m,n}$ 
such that the $(h+k)$-tuple $(a_1,\ldots, a_{h+k-1},a_{h+k}) \in I_{m,n}^{h+k}$ 
satisfies equation~\eqref{quad:DEc}.  
The number of integers that satisfy inequality~\eqref{quad:mn1} is at least
\[
[\beta n] - [\alpha n] > (\beta - \alpha ) n -1
\]
and the  number of integers that satisfy inequality~\eqref{quad:mn2} is at least
\[
[\delta n] - [\gamma n] > (\delta - \gamma ) n - 1.
\]
If
\[
0 < \theta <  (\beta - \alpha )^h (\delta - \gamma )^{k -1} 
\]
then 
\begin{align*}
\Phi^{(c)}_{h,k,m}(n) & \geq ([\beta n] - [\alpha n] )^h ([\delta n] - [\gamma n] )^{k-1} \\
& > ((\beta - \alpha ) n -1)^h ((\delta - \gamma ) n - 1)^{k -1} \\
& = (\beta - \alpha )^h (\delta - \gamma )^{k -1} n^{h+k -1} - O\left( n^{h+k -2}  \right) \\
& > \theta n^{h+k -1} 
\end{align*}
for  all sufficiently large $n$.
This completes the proof.  
\end{proof}

\bl     \label{quad:lemma:fps}  
Let $N$ be a positive integer.  
In the field of formal Laurent series with rational coefficients, if
\[
F_n(z) = \sum_{i=0}^{n-1} z^i = \frac{1 - z^n }{1 - z}
\]
then 
\[
F_n(z)^N
=  \sum_{w=0}^{\infty}  \quad
\sum_{ u=0}^{ \min(N,[w/n]) } 
 (-1)^u \binom{N}{u} \binom{ w - un + N -1}{ N -1} z^w .
\]
\el

\begin{proof}
Recall the binomial theorem
\[
(1-z)^N = \sum_{u=0}^N (-1)^u \binom{N}{u}z^u.
\]
From the formal power series for the $(N-1)$st
derivative of $1/(1-z) = \sum_{v=0}^{\infty}z^v$, we obtain 
\[
\frac{1}{(1-z)^N} = \sum_{v=0}^{\infty} \binom{v+N-1}{N-1} z^v
\]
and so 
\begin{align*}
F_n(z)^N 
& = \frac{   (1-z^n)^N} {(1-z)^N } \\
& =  \left( \sum_{u=0}^{N} (-1)^u \binom{N}{u}z^{un} \right)\left( \sum_{v=0}^{\infty}\binom{ v + N -1}{N -1} z^v   \right) \\
& = \sum_{u=0}^{N} \sum_{v=0}^{\infty}  (-1)^u \binom{N}{u} \binom{ v + N -1}{N -1} z^{un + v} \\
 & =   \sum_{w=0}^{\infty}  \quad \sum_{ \substack{  v\geq 0 \\ 0 \leq u \leq N \\ un+v = w }}(-1)^u \binom{N}{u} \binom{ v + N -1}{N -1} z^w  \\
& =   \sum_{w=0}^{\infty}  \quad \sum_{ u=0}^{ \min(N, [w/n] ) } 
(-1)^u \binom{N}{u} \binom{ w-un + N -1}{N -1} z^w. 
\end{align*}
This completes the proof.  
\end{proof}

\bt             \label{quad:theorem:(k,l,m))poly}
For all $h,k,n \in \N$ and $m,c \in \Z$, 
\[
\Phi^{(c)}_{h,k,m}(n) =  \sum_{ u=0}^{u_0} 
 (-1)^u \binom{h+k}{u} \binom{   (k - u)n - (h - k) m + h + c - 1  }{h + k -1}
\]
where
\begin{align*}
u_0 & = \min\left( h+k,  k +\left[ \frac{ c - k - (h - k) m }{n} \right] \right) \\
& = \begin{cases}
k - 1 & \text{ if $ c-k -  (h - k) m < 0$ and $n \geq  |c-k -  (h - k) m| $} \\
k & \text{ if $ c-k -  (h - k) m  \geq 0$ and $n \geq  c-k - (h - k) m + 1$.} 
\end{cases}
\end{align*}
and 
$\Phi^{(c)}_{h,k,m}(n)$ is a polynomial in $n$ of degree $h + k-1$ 
for sufficiently large $n$.
\et

\begin{proof}
We continue to work in the field of formal Laurent series with rational coefficients.  
Let 
\[
G_{m,n}(z) = \sum_{i=m}^{m+n-1} z^i = z^m F_n(z).
\]
Then
\[
G_{m,n}\left(  \frac{1}{z} \right) = \sum_{i=m}^{m+n-1}  \frac{1}{z^i} 
=   \frac{1}{z^{m+n-1}} \sum_{i=0}^{n-1} z^i 
=   \frac{1}{z^{m+n-1}} F_n(z).
\]
Applying Lemma~\ref{quad:lemma:fps} with $N = h+k$, we have 
\begin{align*}
G_{m,n}(z)^h & G_{m,n}\left(  \frac{1}{z} \right)^{k} 
 = \frac{z^{  hm  }}{z^{ k (m+n-1) }} F_n(z)^{h+k} \\
& = \frac{z^{  hm  }}{z^{ k (m+n-1) }}  \sum_{w=0}^{\infty}  \quad
\sum_{ u=0}^{ \min( h+k, [w/n]) } 
 (-1)^u \binom{h + k}{u} \binom{ w - un + h + k -1}{h + k -1} z^w \\
 & = \sum_{w=0}^{\infty}   \sum_{ u=0}^{ \min( h+k , [w/n] ) } 
 (-1)^u \binom{h + k}{u} \binom{ w - un + h + k -1}{h + k -1} z^{w + (h- k)m - k(n-1)}.
\end{align*}
The coefficient of $z^c$ in this formal Laurent series is $\Phi^{(c)}_{h,k,m}(n)$, 
the number of representations of $c$ 
in the form $a_1+\cdots + a_h - a_{h+1} - \cdots - a_{h+k}$.  
If 
\[
w + (h - k) m - k (n-1) = c
\]
then 
\[
w = kn + c - k - (h-k)m 
\]
and 
\[
w - un + h + k -1 = (k - u)n - (h - k) m + h + c - 1.
\]
Moreover, 
\[
\frac{w}{n} = k + \frac{c - k  - (h-k)m }{n}.
\]
It follows that 
\[
\Phi^{(c)}_{h,k,m}(n) 
=  \sum_{u=0}^{u_0} 
 (-1)^u \binom{h+k}{u} \binom{  (k - u)n - (h - k) m + h + c - 1 }{h + k -1} 
\]
where
\begin{align*}
u_0 & =  \min \left( h+k, \left[  \frac{w}{n}\right] \right) \\
& = \min\left( h+k,  k + \left[  \frac{ c - k - (h - k) m }{n} \right] \right)  \\
& = \begin{cases}
k - 1 & \text{ if $c - k  - (h - k) m  <  0$ and $n \geq |c - k - (h - k) m| $} \\
k & \text{ if $c - k - (h - k) m \geq 0$ and $n \geq c - k - (h - k) m + 1$.} 
\end{cases}
\end{align*}
Each term in this finite series for $\Phi^{(c)}_{h,k,m}(n)$
is a polynomial in $n$ of degree $h+k -1$, 
and so $\Phi^{(c)}_{h,k,m}(n) $  is a polynomial in $n$ of degree 
at most $h+k - 1$.
Lemma~\ref{quad:lemma:polygrowth} implies that $\Phi^{(c)}_{h,k,m}(n)$  
is a polynomial in $n$ of degree 
exactly $h+k - 1$.
This completes the proof.  
\end{proof}

Note that Theorem~\ref{quad:theorem:(k,l,m))poly} is consistent with the 
previous observation that if $h = k$, then 
$\Phi^{(c)}_{h,h, m}(n)= \Phi^{(c)}_{h,h}(n)$ for all $m \in \Z$.

\section{Examples}

\subsection{Additive quadruples and other additive $h$-tuples}

This is the case  $h=k$ and $m=c=0$.  
Applying Theorem~\ref{quad:theorem:(k,l,m))poly}, 
we obtain 
\[
\Psi_h(n) =  \sum_{ u=0}^{h-1} 
 (-1)^u \binom{2h}{u} \binom{   (h - u)n + h - 1  }{2h -1}
\]
for  $n \geq h$.  
Computing with Maple these polynomials for $h = 1,\ldots, 7$, we obtain
\begin{align*}
\Psi_1(n) & = n \\
\Psi_2 (n) & = \frac{2}{3} \,{n}^{3}+\frac{1}{3} \, n \\
\Psi_3(n) & = {\frac {11}{20}} \, {n}^{5}+\frac{1}{4} \, {n}^{3}+\frac{1}{5} \, n \\
\Psi_4 (n) & =    {\frac {151}{315}} \, {n}^{7}+\frac{2}{9}\, {n}^{5}+{\frac {7}{45}} \, {n}^{3}+\frac{1}{7} \, n \\
\Psi_5 (n) & =   {\frac {15619}{36288}} \, {n}^{9}+{\frac {175}{864}} \, {n}^{7}
+{\frac {247}{1728}}\,{n}^{5} +{\frac {1025}{9072}}\,{n}^{3}+ \frac{1}{9} \,n \\
\Psi_6 (n) & =   {\frac {655177}{1663200}}\,{n}^{11}+{\frac {809}{4320}}\,{n}^{9} +{\frac {961}{7200}}\,{n}^{7}+{\frac {3197}{30240}}\,{n}^{5} 
+{\frac {479}{5400}}\,{n}^{3}+\frac{1}{11}\,n    \\
\Psi_7 (n) & =   {\frac {27085381}{74131200}}\,{n}^{13}+{\frac {30233}{172800}}\,{n}^{11}+{\frac {21707}{172800}}\,{n}^{9}+{\frac {51937}{518400}}\,{n}^{7}+{\frac {1813}{21600}}\,{n}^{5}+{\frac {2891}{39600}}\,{n}^{3}+\frac{1}{13} \,n . 
\end{align*}
These polynomials contain only odd powers of $n$.  

The polynomial $\Psi_3(n)$  is in the Online Encyclopedia of Integer Sequences, 
A071816, along with the generating function
\[
\sum_{n=0}^{\infty} \Psi_3(n)   x^n = \frac{x( 1 + 14 x + 36x^2 + 14x^3 + x^4 )}{(1-x)^6}.
\]

\subsection{The case $(h,k) = (3,2)$}
Theorem~\ref{quad:theorem:(k,l,m))poly} produces the polynomial 
\[
\Phi^{(c)}_{3,2,m}(n) =  \sum_{ u=0}^{u_0} 
 (-1)^u \binom{5}{u} \binom{   (2 - u)n -  m + 2 + c  }{4}
\]
where
\begin{align*}
u_0 
& = \begin{cases}
1 & \text{ if $ c \leq  m + 1$ and $n \geq   |c - m - 2| $} \\
2 & \text{ if $ c \geq m + 2$ and $n \geq  c - m - 1$.} 
\end{cases}
\end{align*}
If $m = 0$ and $c \leq 1$, then $u_0 = 1$ and 
\[
\Phi^{(c)}_{3,2}(n)
 =  \binom{  2n  + 2 + c  }{4} - 5 \, \binom{  n  + 2 + c  }{4}
\]
For example, for $m = c = 0$, we obtain 
\[
\Phi^{(0)}_{3,2}(n)  = \binom{2n+2}{4} - 5 \binom{n+2}{4} 
= \frac{11}{24} n^4  + \frac{1}{4} n^3 + \frac{1}{24} n^2  +\frac{1}{4} n.
\]
If $m = 0$ and $c \geq 2$, then $u_0 = 2$ and 
\[
\Phi^{(c)}_{3,2}(n)
 =  \binom{  2n  + 2 + c  }{4} - 5 \, \binom{  n  + 2 + c  }{4} + 10 \, 
 \binom{2+c}{4}.
\]
For  $m = 0$ and $c = 12$, we obtain 
\begin{align*}
\Phi^{(12)}_{3,2}(n) 
& = \binom{2n+14}{4} - 5 \binom{n+14}{4}  + 10 \binom{14}{4} \\
& = \frac {11}{24}{n}^{4}+{\frac {25}{4}}{n}^{3}-{\frac {935}{24}}{n}^{2} 
-{\frac {3875}{4}}n+6006
\end{align*}
 for $n \geq 11$.  
Moreover, $n=11$ is the smallest value of $n$ for which the polynomial gives 
the correct number of solutions of the equation $a_1+a_2+a_3-a_4-a_5 = 12$ with
$(a_1,a_2,a_3,a_4,a_5) \in I_n^5$.

\section{Additive energy and open problems} 

Let $A$ be a nonempty finite subset of an additive abelian group.  
The \emph{$h$-additive energy} of $A$ is
\[
\omega_h(A) = \frac{\Psi_h(A)}{|A|^{2h-1}}.
\]
The \emph{$h$-fold sumset} of $A$ is the set
\[
hA = \{a_1+\cdots + a_h : a_i \in A \text{ for } i=1,\ldots, h\}.
\]
The \emph{$h$-fold Freiman constant} of the set $A$ is 
\[
\kappa_h(A) = \frac{ |hA|}{|A|}.
\]
For $x \in hA$, let $r_{A,h}(x)$ denote the number of $h$-tuples 
$(a_1,\ldots, a_h) \in A^h$ such that 
$a_1+\cdots + a_h = x$.
From the definition of the function $\Psi_h(A)$, we obtain 
\[
\sum_{x\in hA} r_{A,h}(x)^2 = \Psi_h(A).
\]
An application of the Cauchy-Schwarz inequality to the identity 
\[
|A|^h = \sum_{x\in hA} r_{A,h}(x) 
\]
yields 
\[
|A|^{2h} = \left( \sum_{x\in hA} r_{A,h} (x)  \right)^2 
\leq |hA|  \sum_{x\in hA} r_{A,h} (x)^2 
= \kappa_h(A) |A| \Psi_h (A).
\]
Division by $|A|^{2h}$ gives the following  \emph{uncertainty inequality} 
in additive number theory (cf. Nathanson~\cite{nath13c}):
\[
\kappa_h(A) \omega_h(A) \geq 1.
\]

For the interval $I_n$ in the additive group \Z, we have $hI_n  = I_{hn - h+1}$ and so 
\[
\kappa_h(I_n) = \frac{| hI_n |}{ | I_n|} = \frac{hn-h+1}{n} = h - \frac{h-1}{n}.
\]
If $h=2$, then 
\[
\omega_2(I_n) = \frac{\Psi_2(n)}{n^3} = \frac{2}{3}+\frac{1}{3n^2} 
\]
and 
\[
\kappa_2(I_n) \omega_2(I_n) 
=    \left( 2 - \frac{1}{n} \right)  \left( \frac{2}{3}+\frac{1}{3n^2} \right)
= \frac{4}{3} - \frac{ n^2 + (n-1)^2 }{3n^3} > 1
\]
for $n \geq 2$.
If $h=3$, then 
\[
\omega_3(I_n) = \frac{\Psi_3(n)}{n^5} 
= \frac{11}{20}+\frac{1}{4n^2} +\frac{1}{5n^4} 
\]
and 
\begin{align*}
\kappa_3(I_n) \omega_3(I_n) 
& =    \left( 3 - \frac{2}{n} \right)  \left( \frac{11}{20}+\frac{1}{4n^2} +\frac{1}{5n^4} \right) \\
& =  \frac{33}{20} - \frac{22n^4-15n^3+10n^2-12n+8   }{20 n^5}  > 1
\end{align*}
for $n \geq 2$.

Let  $\ell(h)$ denote the leading coefficient of the  polynomial $\Psi_h(n)$.
Then
\[
U(h) = \lim_{n\rightarrow \infty} \kappa_h(I_n) \omega_h(I_n) = h \ell(h).
\]
The following table contains  $\ell(h)$ and $U(h)$ for various values of $h$: \\

\begin{center}
\begin{tabular}{ | c | c | c | c | c | c | c | c  |} \hline 
$h$ & 1&2 & 3 &  4& 5&6 & 7   \\ \hline 
$\ell(h)$ & 1.0000 & 0.66667& 0.55000  & 0.47936 & 0.43042 &0.39393  & 0.36537  \\
\hline
$U(h)$ &  1.0000 & 1.3333   & 1.6500 &1.9174  & 2.1521 & 2.3636 & 2.5576  \\
\hline 
\end{tabular}

\vspace{0.5cm}

\begin{tabular}{ | c | c | c | c | c | c | c | c  |} \hline 
$h$ &  8 & 9 & 10 & 20 &  30 & 40 & 50  \\ \hline 
$\ell(h)$ &  0.34224  & 0.32301& 0.30669 & 0.21769 & 0.17797 & 0.15422 &  0.13799 \\
\hline
$U(h)$ &  2.7379 & 2.9071& 3.0669  &  4.3538 &  5.3391 & 6.1688 & 6.8995 \\
\hline 
\end{tabular}
\end{center}
\vspace{0.5cm}

These calculations suggest several questions.
\benum
\item
Are the sequences $(\ell(h))_{h=1}^{\infty}$ and $(U(h))_{h=1}^{\infty}$ 
monotonically decreasing and increasing, respectively?
If so, what are their limits?

\item
Let $(A_n)_{n=1}^{\infty}$ be a sequence of nonempty finite sets of 
integers such that $ \lim_{n\rightarrow \infty}|A_n| = \infty$.
Is it true that 
\[
\liminf_{n\rightarrow \infty} \kappa(A_n) \omega(A_n) > 1?
\]
If the answer is ``yes,''  then does there exist $\delta > 0$ such that
\[
\liminf_{n\rightarrow \infty} \kappa(A_n) \omega(A_n) \geq 1+ \delta
\]
for all such sequences $(A_n)_{n=1}^{\infty}$?

\item
Consider the linear form
\[
L(x_1,\ldots, x_h) = \lambda_1x_1 + \cdots + \lambda_h x_h
\]
with integer coefficients $\lambda_1, \ldots, \lambda_h$.  
Determine the polynomiality of the number of solutions of $L(x_1,\ldots, x_h) = c$
with $(x_1,\ldots, x_h) \in I_n^h$ as $n \rightarrow \infty$.

\item
Additive $(h,k)$-tuples in finite fields occur often in number theory.  
For example, a solution of the equation 
$a_1 + a_2 + a_3 + a_ 4 = a_5 + a_6 + a_7 + a_ 8$
is called an \emph{additive octuple}.   
Bateman and Katz~\cite{bate-katz11} have considered additive octuples 
 in finite fields.
 Morier-Genoud and Ovsienko~\cite{mori-ovsi13} consider additive quadruples 
 in the vector space $\F_2^n$ in connection with Hurwitz sums of squares identities.  
 Green~\cite{gree05} and Green and Tao~\cite{gree-tao09} also consider 
 additive quadruples in finite fields.
It would of interest to determine the polynomiality of the number of additive 
quadruples, additive octuples, and other additive $(h,k)$-relations 
in the vector space $\F_q^n$ over the finite field $\F_q$ as 
the dimension $n \rightarrow \infty$.

\eenum

\def\cprime{$'$} \def\cprime{$'$} \def\cprime{$'$}
\providecommand{\bysame}{\leavevmode\hbox to3em{\hrulefill}\thinspace}
\providecommand{\MR}{\relax\ifhmode\unskip\space\fi MR }
\providecommand{\MRhref}[2]{%
  \href{http://www.ams.org/mathscinet-getitem?mr=#1}{#2}
}
\providecommand{\href}[2]{#2}

\end{document}